\documentclass{article}

\usepackage[margin=1in]{geometry}
\usepackage{amsmath, amssymb, tikz-cd, centernot, mathtools, csquotes,
            mathrsfs, amsthm, xparse, microtype}
\usepackage[utf8]{inputenc}
\usepackage[T1]{fontenc}
\usepackage[backend=biber, style=numeric, backref=true]{biblatex}
\usepackage[nameinlink]{cleveref}
\usepackage{IEEEtrantools}

\newcounter{smarttheorem:smarttheoremcounter}
\newcommand{\smarttheoremlabel}{%
  smarttheorem:\arabic{smarttheorem:smarttheoremcounter}%
}
\NewDocumentCommand{\newsmarttheorem}{m o m}{
  \IfNoValueTF{#2}{\newcounter{#1}}{}
  \newtheorem{hidden#1}[\IfValueTF{#2}{#2}{#1}]{#3}

  \NewDocumentEnvironment{#1}{o}{
    \IfNoValueTF{##1}{\begin{hidden#1}}{\begin{hidden#1}[{##1}]}
      \stepcounter{smarttheorem:smarttheoremcounter}
      \label{\smarttheoremlabel}
  }{
      \end{hidden#1}
  }
}

\let \vqedsymbol \qedsymbol

\NewDocumentEnvironment{pf}{o o}{
  \IfValueTF{#2}{
    \edef \qedsymbollabel {#2}
    \renewcommand \qedsymbol {\vqedsymbol \, \qedsymbollabel}
  }{
    \renewcommand \qedsymbol \vqedsymbol
  }
  \IfValueTF{#1}{\begin{proof}[#1]}{\begin{proof}}
}{
  \end{proof}
}

\NewDocumentEnvironment{rpf}{O{\proofname} O{\smarttheoremlabel}}{
  \begin{pf}[#1][\noexpand\Cref{#2}]
}{
  \end{pf}
}

\NewDocumentEnvironment{lrpf}{m}{
  \begin{rpf}[Proof of \Cref{#1}][#1]
}{
  \end{rpf}
}

\newcommand{\plr}{\item[\((\implies)\)]}
\newcommand{\prl}{\item[\((\impliedby)\)]}

\NewDocumentEnvironment{ea}{O{rCl}}{
  \begin{IEEEeqnarray*}{#1}
}{
  \end{IEEEeqnarray*}
  \ignorespacesafterend
}

\NewDocumentEnvironment{tcd}{s}{
  \IfBooleanTF{#1}{\begin{equation}}{\begin{equation*}}
    \begin{tikzcd}
}{
    \end{tikzcd}
  \IfBooleanTF{#1}{\end{equation}}{\end{equation*}}
  \ignorespacesafterend
}

\theoremstyle{definition}

\newcounter{case}
\newenvironment{caselist}{
  \setcounter{case}{0}\begin{description}
}{
  \end{description}
}

\newcommand \case {\setcounter{case}{\value{case}+1}\item[Case \thecase.]}

\numberwithin{theoremcounter}{section}
\theoremstyle{plain}
\newsmarttheorem{theorem}[theoremcounter]{Theorem}
\newsmarttheorem{lemma}[theoremcounter]{Lemma}
\newsmarttheorem{proposition}[theoremcounter]{Proposition}
\newsmarttheorem{claim}[theoremcounter]{Claim}
\newsmarttheorem{subclaim}[theoremcounter]{Subclaim}
\newsmarttheorem{corollary}[theoremcounter]{Corollary}
\newsmarttheorem{fact}[theoremcounter]{Fact}
\newsmarttheorem{todo}{TODO}

\theoremstyle{definition}
\newsmarttheorem{definition}[theoremcounter]{Definition}
\newsmarttheorem{notation}[theoremcounter]{Notation}

\theoremstyle{remark}
\newsmarttheorem{question}[theoremcounter]{Question}
\newsmarttheorem{remark}[theoremcounter]{Remark}
\newsmarttheorem{aside}[theoremcounter]{Aside}
\newsmarttheorem{conjecture}[theoremcounter]{Conjecture}
\newsmarttheorem{exercise}[theoremcounter]{Exercise}
\newsmarttheorem{example}[theoremcounter]{Example}

\renewcommand \epsilon \varepsilon
\renewcommand \phi \varphi

\def\Ind#1#2{#1\setbox0=\hbox{$#1x$}\kern\wd0\hbox to 0pt{\hss$#1\mid$\hss}
\lower.9\ht0\hbox to 0pt{\hss$#1\smile$\hss}\kern\wd0}

\def\notind#1#2{#1\setbox0=\hbox{$#1x$}\kern\wd0
\hbox to 0pt{\mathchardef\nn=12854\hss$#1\nn$\kern1.4\wd0\hss}
\hbox to
0pt{\hss$#1\mid$\hss}\lower.9\ht0 \hbox to 0pt{\hss$#1\smile$\hss}\kern\wd0}

\DeclarePairedDelimiter \abs \lvert \rvert

\DeclarePairedDelimiter \pars ( )

\DeclarePairedDelimiterX \set [1] \{ \} {\, #1 \,}

\DeclareMathOperator \thy {Th}
\DeclareMathOperator \tp {tp}
\DeclareMathOperator \rkm {RM}

\DeclareMathOperator \mult {mult}

\newcommand \R {\mathrm{R}}

\title{Remarks on recognizable subsets and local rank}
\author{Christopher Hawthorne\thanks{%
  This work was partially supported by an NSERC CGS-M, an OGS, and an NSERC
  PGS-D.}
}
\date{May 2019}

\begin{document}
\maketitle
\begin{abstract}
  Given a monoid \((M,\epsilon,\cdot  )\) it is shown that a subset
  \(A\subseteq M\) is recognizable in the sense of automata theory if and only
  if the \(\phi \)-rank of \(x=x\) is zero in the first-order
  theory \(\thy(M,\epsilon ,\cdot ,A)\), where \(\phi (x;u)\) is the formula
  \(xu\in A\).  In the case where \(M\) is a finitely generated free monoid on
  a finite alphabet \(\Sigma \), this gives a model-theoretic characterization
  of the regular languages over \(\Sigma \). If \(A\) is a regular language
  over \(\Sigma \) then the \(\phi \)-multiplicity of \(x=x\) is the state
  complexity of \(A\). Similar results holds for \(\phi' (x;u,v)\) given by
  \(uxv\in A\), with the \(\phi' \)-multiplicity now equal to the size of the
  syntactic monoid of \(A\). 
\end{abstract}

\section{Introduction}
There is known to be interplay between logic and the study of regular
languages: for example, if one views words as structures in a suitable
signature, one can characterize the regular languages as the sets of words that
are axiomatizable in monadic second-order logic. (See \cite[Chapter
IX.3]{pin16}.) In this paper we present a
different logical characterization of regularity. Given
a finite alphabet \(\Sigma \) we characterize the regular languages over
\(\Sigma \) by studying the free monoid on \(\Sigma \) as a first-order
structure \((\Sigma ^*,\epsilon ,\cdot )\); in particular we study local
stability theoretic ranks in \(\thy(\Sigma ^*,\epsilon ,\cdot )\). It follows
from the work of Quine that the regular languages over \(\Sigma \) are
definable in \((\Sigma ^*,\epsilon ,\cdot )\). We show that a definable subset
\(A\subseteq \Sigma ^*\) is regular if and only if the \(\phi \)-rank of
\(x=x\) is zero, where \(\phi (x;u)\) is the formula \(xu\in A\); this is
\Cref{cor:regchar} below. 
In particular, when \(A\) is regular \(\phi \) is stable. Moreover, when \(A\)
is regular the \(\phi \)-multiplicity of \(x=x\) agrees
with the state complexity of \(A\).  In fact, our characterization applies to
arbitrary (not necessarily free) monoids, where regularity is replaced by
\emph{recognizability}; see \Cref{thm:RecogTypes} below.

Our approach here is a familiar one in local stable group theory:
one studies a subset \(A\) of the group in question by studying the 
formula \(xy\in A\). See for example \cite{conant18}.

We also take this opportunity to point out (using the known structure of the
atomically definable subsets of \(\Sigma ^*\)) that every quantifier-free definable
subset of \(\Sigma ^*\) is a star-free regular language. Finally, we
investigate some examples demonstrating that the connections between formal
languages and model theory made here cannot be taken much further.

I warmly thank my advisor, Rahim Moosa, for excellent guidance, thorough
editing, and many helpful discussions. The results here appeared as part of my
Master's essay at the University of Waterloo.

\section{A characterization of recognizability}
The following is a basic notion in automata theory. We suggest 
\cite[Section III.12]{eilenberg74} and \cite[Section
IV.2]{pin16} for further details.
\begin{definition}
  Suppose \(M\) is a monoid; suppose \(A \subseteq M\). We say \(A\) is a
  \emph{recognizable subset of \(M\)} if there is some finite monoid \(F\) and
  some homomorphism of monoids \(\alpha \colon M \to F\) such that \(A\) is a
  union of fibres of~\(\alpha\); i.e.\ \(A = \alpha^{-1}(\alpha(A))\). 
  
  A standard fact from automata theory (see \cite[Proposition
  III.12.1]{eilenberg74}) is that \(A\) is recognizable if and
  only if \(\theta _A\) has finitely many equivalence classes, where \(\theta
  _A\) is the equivalence relation on \(M\) given by \((a,b)\in\theta _A\) if
  and only if for all \(c\in M\) we have \(ac\in A\iff bc\in A\).
\end{definition}
Our goal in this paper is to give a model-theoretic characterization of
recognizable subsets of a monoid. The characterization will be in terms of
local model-theoretic ranks, which we now briefly recall; see also
\cite[Exercise 8.2.10]{tent12}.

Suppose we are given a signature \(L\), a complete theory \(T\),
and an \(L\)-formula \(\phi(\overline{x}; \overline{u})\).  We let
\(\mathfrak{C}\) denote a fixed sufficiently saturated model of \(T\).
Suppose \(A\) is a parameter set.
A \emph{\(\phi\)-formula over \(A\)} is a 
Boolean combination of formulas of the form \(\phi
(\overline{x};\overline{a})\) for some tuple \(\overline{a}\) from \(A\).
A \emph{complete \(\phi\)-type over \(A\)} is a maximally consistent set
of \(\phi\)-formulas over \(A\). We use \(S_\phi(A)\) to denote the set
of complete \(\phi\)-types over \(A\).
Suppose \(\overline{a}\) is a tuple from \(\mathfrak{C}\) of the same
arity as \(\overline{x}\). We define the \emph{\(\phi\)-type of
\(\overline{a}\) over \(A\)} to be
$\tp_\phi(\overline{a}/A) = \set{\psi(\overline{x}) : \psi \text{ is a
  \(\phi\)-formula over } A, \mathfrak{C} \models \psi(\overline{a})}
$
\begin{definition}
  Suppose \(\psi(\overline{x})\) is a \(\phi\)-formula. The
  \emph{\(\phi \)-rank} \(R_\phi
  (\psi )\) is defined as follows:
  \begin{itemize}
    \item
      \(R_\phi (\psi )\ge 0\) if \(\psi \) is consistent.
    \item
      If \(\beta \) is a limit ordinal, then \(R_\phi (\psi )\ge\beta \) if
      \(R_\phi (\psi )\ge \alpha \) for all \(\alpha <\beta \).
    \item
      If \(\alpha \) is an ordinal, then \(R_\phi(\psi) \ge \alpha+1\) if there
      are \(\phi\)-formulas \((\psi_i(\overline{x}) : i < \omega)\) satisfying:
      \begin{enumerate}
        \item
          The \(\psi_i\) are pairwise inconsistent.
        \item
          \(\mathfrak{C} \models \forall \overline{x} (\psi_i(\overline{x})
          \rightarrow \psi(\overline{x}))\) for all \(i < \omega\).
        \item
          \(R_\phi(\psi_i) \ge \alpha\) for all \(i < \omega\).
      \end{enumerate}
  \end{itemize}
  If \(R_\phi (\psi )\ge\alpha \) for every ordinal \(\alpha \), we write
  \(R_\phi (\psi )=\infty \). This notion is sometimes referred to by the name
  ``local rank''.
\end{definition}
\begin{definition}
  Suppose \(\psi(\overline{x})\) is a \(\phi\)-formula with \(R_\phi(\psi) <
  \infty\). The \emph{\(\phi \)-multiplicity} \(\mult_\phi(\psi)\) is defined
  to be the largest \(n < \omega\) such that there are \(\phi \)-formulas
  \((\psi_i : i < n)\) satisfying:
  \begin{enumerate}
    \item
      The \(\psi_i\) are pairwise inconsistent.
    \item
      \(\mathfrak{C} \models \forall \overline{x} (\psi_i(\overline{x})
      \rightarrow \psi(\overline{x}))\) for all \(i < n\).
    \item
      \(R_\phi(\psi_i) \ge R_\phi(\psi)\) for all \(i < n\).
  \end{enumerate}
\end{definition}
We work in \(L_P = \set{1,\cdot ,P}\), the signature of monoids
expanded by a unary predicate \(P\).  Given a monoid \(M\) and some \(A
\subseteq M\), we can view \((M,A)\) as an \(L_P\)-structure by interpreting
\(P\) as the set \(A\).

We now present our main theorem:
\begin{theorem}
  \label{thm:RecogTypes}
  Suppose \(M\) is a monoid and \(A\subseteq M\); let \(\phi (x;u)\) be the
  \(L_P\)-formula \(P(xu)\). Then \(A\) is recognizable if and only if
  \(R_\phi(x = x)\) is zero in \(\thy(M,A)\).  Moreover, when this is the case,
  \(\mult_\phi(x = x)\) is the number of equivalence classes of \(\theta_A\).
\end{theorem}
Recall that globally \(\rkm(x = x)\) is zero if and only if the
domain is finite; hence this theorem can be understood as saying that \(A\) is
recognizable if and only if the monoid is ``finite with respect to the formula
\(xu\in A\)''.
\begin{rpf}
  Recall that \(A\) is recognizable if and only if \(\theta _A\) has finitely
  many equivalence classes.  Observe that for \(a,b \in M\) we have \((a,b) \in
  \theta_A\) if and only if \((M,A) \models \phi(a; c) \leftrightarrow \phi(b;
  c) \) for all \(c \in M\).
  This is in turn equivalent to requiring that
  \(\tp_\phi(a/M) = \tp_\phi(b/M)\). Thus \(A\) is recognizable
  if and only if
  only finitely many complete \(\phi \)-types over \(M\) are realized in \(M\);
  it therefore suffices to check that this
  is equivalent to \(R_\phi (x=x)\) being zero.
  \begin{description}
    \plr
      Note for \(a,b\in M\) that
      \(\tp_\phi (a/M) = \tp_\phi (b/M)\) if and only if \((M,A)\models \forall
      z(\phi (a;z)\leftrightarrow \phi (b;z))\).
      So \(\thy(M,A)\) witnesses that only finitely many elements of \(S_\phi
      (M)\) are realized in \(M\),
      and hence the same holds true in all models. It follows that \(S_\phi
      (\mathfrak{C})\) is finite; else we could find an elementary extension
      realizing infinitely complete \(\phi \)-types. But if there are only
      finitely many complete global \(\phi \)-types then there cannot be
      infinitely many consistent and pairwise inconsistent \(\phi \)-formulas;
      it follows that \(R_\phi (x=x)\) is zero.
    \prl
      Suppose there are infinitely many elements of \(S_\phi (M)\) that are
      realized in \(M\).  We will show that \(R_\phi (x=x)\ge 1\).

      We iteratively define \(\phi \)-formulas \((\psi _n: n<\omega) \)
      satisfying the following:
      \begin{itemize}
        \item
          For each \(n < \omega\) we have \(\set{\psi _i : i < n}\)
          extends to infinitely many \(\phi\)-types over \(M\).
        \item
          For each \(n < \omega\) we have that
          \(\set{\psi _i : i < n} \cup \set{\neg \psi_n}\) is consistent.
      \end{itemize}
      Having defined \(\psi_i\) for \(i < n\), we observe that since
      \(\set{\psi _i : i < n}\) extends to infinitely many
      \(\phi\)-types over \(M\), in particular it does not isolate a type, and
      there is some \(\phi \)-formula \(\psi_n \) with
      \(\set{\psi _i : i < n} \cup \set{\psi _n}\) and
      \(\set{\psi _i : i < n} \cup \set{\neg \psi _n}\) both consistent.
      Furthermore, at least one of
      \(\set{\psi _i : i < n} \cup \set{\psi _n}\) and
      \(\set{\psi _i : i < n} \cup \set{\neg \psi _n}\)
      must extend to
      infinitely many \(\phi\)-types over \(M\); without loss of generality, we
      may assume it is the former. (Note that this construction works perfectly
      well in the case \(n = 0\).)

      For \(n<\omega \) we now define
      \[\delta _n = \neg\psi _n\wedge \bigwedge _{i < n}\psi _i
      \]
      Then by construction the
      \(\delta_n\) are consistent and pairwise inconsistent; so \(R_\phi(x = x)
      \ge 1\), as desired.
  \end{description}
  For the ``moreover'' clause, recall from the beginning of the proof that the
  number of complete \(\phi \)-types over \(M\) that are realized in \(M\) is
  equal to the number of equivalence classes of \(\theta _A\); call this number
  \(m\). A careful examination of the left-to-right proof above yields that
  \(m=\abs{S_\phi (\mathfrak{C})}\ge \mult_\phi (x=x)\); it remains to show
  that \(m\le \mult_\phi (x=x)\). Write \(S_\phi (\mathfrak{C}) = \set{p_i : i
  < m}\). Since \(S_\phi (\mathbb{C})\) is finite, each \(p_i\) is isolated by
  some \(\phi \)-formula \(\psi _i\). Then the \((\psi _i:i<m)\) are consistent
  and pairwise inconsistent, and thus \(\mult_\phi (x=x) \ge m\).
\end{rpf}
\begin{corollary}
  If \(A\) is recognizable then \(\phi\) is stable.
\end{corollary}
\begin{rpf}
  From the proof above we get that \(S_\phi (\mathfrak{C})\) is finite if \(A\)
  is recognizable; it follows that \(\phi \) is stable.
\end{rpf}
\section{The case of regular languages}
A case of particular interest is when \(M = (\Sigma ^*,\epsilon ,\cdot )\) is
the free monoid on some finite set \(\Sigma \). In this case, the notion of
a recognizable set coincides with that of a \emph{regular language}:
\begin{definition}
  The class of \emph{regular languages} over a finite set \(\Sigma \) is the
  smallest class of subsets of \(\Sigma ^*\) containing the finite sets and
  closed under union, element-wise concatenation, and Kleene star. (The Kleene
  star of \(A\subseteq \Sigma ^*\) is \(A^* = \set{a_0\cdot \cdots \cdot
  a_{n-1} : n<\omega , a_0,\ldots ,a_{n-1}\in A}\).)
\end{definition}
The regular languages can also be characterized as the subsets of \(\Sigma ^*\)
that can be recognized by a finite automaton on the alphabet \(\Sigma \). See
\cite{yu97} for more details on regular languages.
\begin{corollary}
  \label{cor:regchar}
  Suppose \(\Sigma \) is a finite set and \(A\subseteq \Sigma ^*\) is definable
  in \((\Sigma ^*,\epsilon ,\cdot )\). Let \(\phi (x; u)\) be
  the formula in the language of monoids (with parameters in \(\Sigma ^*\))
  expressing that \(xu\in A\). 
  Then \(A\) is a regular language over \(\Sigma \) if and only if
  \(R_{\phi }(x= x)\) is zero.
\end{corollary}
\begin{rpf}
  Working in the \(L_P\)-theory \(\thy(\Sigma ^*,\epsilon ,\cdot ,A)\) we have
  that and \(\phi \) is equivalent to \(P(xu)\).  Since local ranks don't
  change under definitional expansions of the theory, the result follows from
  \Cref{thm:RecogTypes} and the fact that the recognizable subsets of \(\Sigma
  ^*\) are precisely the regular languages over \(\Sigma \); see \cite[Theorem
  4.7]{yu97}.
\end{rpf}
This is at first glance only a statement about the regular languages that also
happen to be definable; but in fact all regular languages are definable.
In \cite{quine46} Quine presents an interpretation of
\((\mathbb{N},+,\cdot )\) in \((\Sigma ^*,\epsilon ,\cdot )\) for any
\(2\le\abs{\Sigma }<\omega \) induced by a computable bijection
\(\Phi \colon \mathbb{N}\to\Sigma ^*\). Given a regular \(A\subseteq \Sigma
^*\), we first note that membership in \(A\) is decidable; this follows easily
from the automata-theoretic characterization of regularity.
Hence \(\Phi ^{-1}(B)\subseteq \mathbb{N}\) is also decidable. But every
decidable subset of \(\mathbb{N}\) is definable in \((\mathbb{N},+,\cdot )\)
(see for example \cite[Theorem 6.12]{flum94}); so \(\Phi ^{-1}(A)\) is
definable in \((\mathbb{N},+,\cdot )\). So \(A = \Phi (\Phi ^{-1}(A))\)
is definable in \((\Sigma ^*,\epsilon ,\cdot)\).

So we have completely characterized the regular languages over \(\Sigma \).
Moreover, when \(A\subseteq \Sigma ^*\) is a regular language, we can
meaningfully interpret the multiplicity.  Indeed, it follows from the
\enquote{moreover} clause of \Cref{thm:RecogTypes} that \(\mult_\phi (x=x)\) is
the \emph{state complexity} of \(A\): the number of states in the smallest
automaton recognizing \(A\). Details about the state complexity and its
connection to \(\theta _A\) can be found in \cite[Sections
3.9-3.11]{shallit08}.
\begin{remark}
  There is another equivalence relation whose classes characterize
  recognizability: the \emph{syntactic congruence}, given by \((a, b)\in\theta
  _A'\) if and only if for all \(c_1,c_2\) we have \(c_1ac_2\in A\iff
  c_1bc_2\in A\).  \Cref{thm:RecogTypes} remains true if we replace \(\theta
  _A\) with \(\theta _A'\) and \(\phi \) with \(\phi' (x; u,v) := P(uxv)\); so
  we can also characterize recognizability by whether \(R_{\phi '} (x=x)\) is
  zero. We can again meaningfully interpret the multiplicity:
  \(\mult_{\phi' }(x=x)\) will be the size of the \emph{syntactic
  monoid} of \(A\), which is the smallest monoid witnessing the recognizability
  of \(A\). Details about syntactic monoids can be found in \cite[Section
  III.10]{eilenberg74} or \cite[Section IV.4]{pin16}.
\end{remark}
Having characterized the regular languages among the definable subsets of
\(\Sigma ^*\), one might ask which quantifier-free definable subsets of
\(\Sigma ^*\) are regular; it turns out that all of them are, and they are
furthermore star-free.
\begin{definition}
  The class of \emph{star-free languages} over an alphabet \(\Sigma \) is the
  smallest class of subsets of \(\Sigma ^*\) that contains the finite sets and
  is closed under Boolean combinations and element-wise concatenation.
  (Note that every star-free language is regular.)
\end{definition}
The following is a consequence of the known structure of the atomically
definable subsets of finitely generated free monoids:
\begin{proposition}
  \label{prop:qfstarfree}
  Every quantifier-free definable subset of \((\Sigma^*,\epsilon ,\cdot )\) is
  star-free, and in particular is regular.
\end{proposition}
\begin{rpf}
  \cite[Theorem 3]{dabrowski2011} tells us that if \(A\subseteq \Sigma ^*\) is
  atomically definable then \(A\) is either finite, all of \(\Sigma ^*\), or is
  the union of a finite set and \(a^*a_0\) for some \(a\in\Sigma ^*\)
  primitive and some proper prefix \(a_0\) of \(a\). (We say \(a\) is
  \emph{primitive} if it can't be written in the form \(b^n\) for some
  \(b\in\Sigma ^*\) and \(n>1\).) Since the class of star-free languages is
  closed under Boolean combinations, it then suffices to check that each of
  these is star-free; since finite sets and \(\Sigma ^*=\emptyset ^c\) are
  star-free, it remains to check that \(a^*a_0\) is star-free if \(a\) is
  primitive. But \cite[Theorem 2.5]{mateescu97} asserts that \(a^*\) is
  star-free if \(a\) is primitive, so since star-free languages are closed
  under concatenation we get that \(a^*a_0\) is star-free, as desired.
\end{rpf}

\section{Some counterexamples}
We finish by answering negatively some questions that arise naturally from the
above considerations.

Having just shown that every quantifier-free definable subset of \(\Sigma ^*\)
is star-free, one might ask whether the converse holds; it does not.
\begin{example}
  Let \(A = 0^*1^* \subseteq \set{0,1}^*\).  Then \(A = (\emptyset
  ^c10\emptyset ^c)^c\) is star-free but not quantifier-free definable.
\end{example}
\begin{rpf}
  Given \(B\subseteq \Sigma ^*\) we let \(n_k(B) = \abs{B\cap \Sigma ^k}\);
  note that \(n_k(A) = k+1\) for all \(k\).  Suppose for contradiction that
  \(A\) is defined by a quantifier-free formula in disjunctive normal form
  \[\bigvee _{i<p}\bigwedge _{j<q_i} \phi _{ij}(x)
  \]
  with each \(\phi _{ij}\) atomic or negated atomic. The description of the
  atomically definable sets mentioned in the proof of \Cref{prop:qfstarfree}
  tells us that for sufficiently
  large \(k\) we have \(n_k(\phi _{ij}(\Sigma ^*))\) is either \(\le 1\) (if
  atomic) or \(\ge 2^k-1\) (if negated atomic). (We may assume that none of the
  \(\phi _{ij}\) is a tautology.) We have two cases:
  \begin{caselist}
    \case
      Suppose for some \(i_0\) we had that each \(\phi _{i_0j}\) was negated
      atomic.
      Then for sufficiently large \(k\) we would have
      \[k+1 = n_k(A) = \abs{A\cap \Sigma^k}  \ge \abs*{\bigcap _{j<q_{i_0}} \phi
        _{i_0j}(\Sigma ^*)\cap \Sigma ^k} \ge 2^k-q_{i_0}
      \]
      a contradiction.
    \case
      Suppose for each \(i\) there is \(j_i\) such that \(\phi _{ij_i}\) is
      atomic. Then for sufficiently large \(k\) we get for each \(i\) that
      \[n_k\pars*{\bigcap _{j<q_i} \phi _{ij}(\Sigma ^*)} \le n_k(\phi
        _{ij_i}(\Sigma ^*)) \le 1
      \]
      so \(k+1=n_k(A)\le m\), a contradiction.
  \end{caselist}
  So no such formula exists, and \(A\) isn't quantifier-free definable.
\end{rpf}
We noted earlier that every recognizable subset of a finitely generated free
monoid is definable in the monoid structure; this doesn't generalize to
arbitrary monoids.
\begin{example}
  \label{eg:recnotdef}
  Let
  \[M=\bigoplus _{i<\omega }(\mathbb{Z}/2\mathbb{Z},+)
  \]
  and let \(A\) be the set of elements of \(M\) with an even number of non-zero
  entries. Then \(A\) is recognizable but not definable in \((M,0,+)\).
\end{example}
\begin{rpf}
  The monoid homomorphism \(M\to\mathbb{Z}/2\mathbb{Z}\) given by
  \[(n_i:i<\omega )\mapsto \sum_{i<\omega }n_i
  \]
  recognizes \(A\).  (Note that this sum is well-defined since there can only
  be finitely many non-zero entries.)
  
  It follows from \cite[Theorem 1.2]{cherlin82} that \(\thy(M)\) admits
  quantifier elimination (this is also not hard to verify by hand); one can
  also check that the only quantifier-free definable subsets of \(M\) are
  finite or cofinite. Since \(A\) is neither, we get that \(A\) is not
  definable.
\end{rpf}

It follows from the proof of \Cref{thm:RecogTypes} that if \(A\) is
a recognizable subset of a monoid \(M\) then every \(\phi\)-type over \(M\) is
realized in \(M\). The following example shows that the converse is false. 
\begin{example}
  Let \(M = \set{0,1,2}^*\).  Given \(a \in M\) and \(i \in \set{0,1,2}\), we
  let \(n_i(a)\) denote the number of occurrences of \(i\) in \(a\). Now, let
  \(A = \set{a \in \set{0,1}^* : n_0(a) = n_1(a)} \subseteq M\); let \(T =
  \thy(M,A)\). Then \((0^n : n < \omega)\) are
  pairwise unrelated by \(\theta_A\), so \(\theta _A\) has infinitely many
  equivalence classes and thus \(A\) is not recognizable (again see
  \cite[Theorem III.12.1]{eilenberg74}); nonetheless every
  complete \(\phi\)-type over \(M\) is realized in \(M\).
\end{example}
\begin{rpf}
  For \(a \in M\), let \(\Delta(a) = n_1(a) - n_0(a)\). Note that for any
  \(a,c \in M\) we have \((M,A) \models \phi(a; c)\) if and
  only if \(a,c
  \in \set{0,1}^*\) and \(\Delta(a) + \Delta(c) = 0\). In
  particular, if \(c\in \set{0,1}^*\) then
  \(\phi(x;c) \wedge \neg\phi(x; d)\) is
  inconsistent if
  \(d \in \set{0,1}^*\) and \(\Delta(c) = \Delta(d)\), and \(\phi
  (x;c)\wedge \phi (x;d)\) is inconsistent if
  \(d\notin \set{0,1}^*\) or \(\Delta (c)\ne\Delta (d)\). Hence
  \(\phi(x; c)\) isolates a
  \(\phi\)-type over \(M\), namely the one determined by
  \[\set{\phi (x;d) : d \in \set{0,1}^*, \Delta (c) = \Delta (d)}\cup
    \set{\neg\phi (x;d) : d\notin \set{0,1}^*\text{ or } \Delta (c)\ne \Delta
    (d)}
  \]
  Hence if \(p(x) \in S_\phi(M)\) contains \(\phi(x; c)\) for
  some \(c \in \set{0,1}^*\), then \(p\) is isolated, and thus realized
  in \(M\).
  Furthermore, if \(c\notin \set{0,1}^*\)
  then \(\phi(x; c)\) is
  inconsistent; hence if \(p(x) \in S_{\phi}(M)\) contains no formula of
  the form \(\phi(x;c)\) for some \(c \in \set{0,1}^*\), then
  \(p\) is the complete \(\phi\)-type over \(M\) determined by
  \(\set{\neg \phi(x; c) : c \in M} \), and \(p\) is realized by \(2 \in M\).
  So every complete \(\phi\)-type over \(M\) is realized in \(M\).
\end{rpf}
The next two examples suggest that while local rank zero coincides with
regularity, it does not seem like higher local ranks correspond to higher
levels in the formal languages hierarchy. We first exhibit a
deterministic context-free language (that is, a very tame non-regular language)
that gives rise to infinite \(\phi \)-rank. (A good reference for
deterministic context-free languages is \cite[Section 4.7]{shallit08}.) In the
other direction, in \Cref{eg:lowranknotcfl} below, we exhibit a language that
is not context-free (so relatively wild) but with respect to which the universe
is of \(\phi \)-rank and multiplicity \(1\).
\begin{example}
  Let \(M=\set{0,1,2}^*\). As a notational convenience, for \((n_i:i<k)\)
  a tuple of natural numbers, we let \((n_i:i<k)_*\) denote
  \(0^{n_0}10^{n_1}1\cdots 10^{n_{k-1}}\); observe that every element of
  \(\set{0,1} ^*\) has a unique representation in this form. Let
  \[A=\set{(n_i:i<k)_*2(m_i:i<\ell )_*:\ell \le k,m_i = n_{k-i-1} + 1 \text{
    for some } i<\ell }
  \]
  It is straightforward to show that \(A\) is a deterministic context-free
  language; nonetheless \(R_{\phi }(x=x)\) is infinite.
\end{example}
\begin{rpf}
  The important property of \(A\) is that if \(i,r\in\mathbb{N}\) and \(a \in
  M\), then \((M,A)\models \phi (a;
  21^j0^r)\) if and only if \(a=(n_i:i<k)_*\) where \(j < k\) and \(n_{k-j-1}
  +1= r\); thus thinking of the elements of \(\set{0,1}^*\) as finite sequences
  of naturals, we see that the components can be extracted via \(\phi
  \)-formulas. This comparison with finite sequences suggests that our
  universe is in some sense not finite dimensional, lending intuitive support
  to the claim that \(R_\phi (x=x)\) is infinite.

  For \(\sigma = (\sigma _i : i < \abs\sigma  )\) a tuple of natural numbers,
  define
  \[\psi _{\sigma }(x) = \bigwedge _{i<\abs\sigma }\phi (x; 21^i0^{\sigma
    _i+1})
  \]
  Thus for non-empty \(\sigma \) the realizations of \(\psi _{\sigma }(x)\) are
  exactly \((n_i:i<k)_*\) where \((n_{k-1},\ldots ,n_{k-\abs\sigma -1})=\sigma
  \). In particular, the \(\psi _\sigma \) are \(\phi \)-formulas
  satisfying:
  \begin{enumerate}
    \item
      \(\psi _{\sigma \hat{} n_1}\) and \(\psi _{\sigma \hat{} n_2}\) are
      inconsistent for all \(\sigma \) and \(n_1\ne n_2\).
    \item
      \((M,A)\models \forall x(\psi _{\sigma \hat{} n}(x)\rightarrow \psi
      _\sigma (x))\) for all \(\sigma\) and all \(n\).
    \item
      \(\psi _\sigma \) is consistent for all \(\sigma \).
  \end{enumerate}
  At this point a simple transfinite induction yields that given such a tree of
  \(\phi \)-formulas, the \(\phi \)-rank of any formula in the tree
  is \(\infty \); hence in particular \(R_\phi (x=x)\) is infinite.
\end{rpf}
\begin{example}
  \label{eg:lowranknotcfl}
  Consider \(M=\set{0,1,2}^*\); let \(A=\set{0^n1^n2^n : n<\omega }\). Then
  \(A\) is not context-free; nonetheless \(R_\phi (x=x)=\mult_\phi (x=x)=1\).
\end{example}
\begin{rpf}
  That \(A\) is not context-free is \cite[Chapter 4, Exercise 1]{shallit08}.
  Note that every consistent \(\phi (x; c)\) with
  \(c\in M\) extends to at most two realized \(\phi \)-types over
  \(M\); this is then witnessed by \(\thy(M,A)\). Thus for every
  \(\mathcal{N}\models \thy(M,A)\) we get for all \(\gamma \in
  N\) that \(\phi (x; \gamma )\)
  extends to at most two realized \(\phi \)-types over \(N\).
  Thus for all \(\gamma \in \mathfrak{C}\) we get that
  \(\phi (x;\gamma )\) extends to at most two \(\phi \)-types (realized or
  otherwise) over \(\mathfrak{C}\) (else there would
  be some elementary extension \(\mathcal{N}\) of \(\mathfrak{C}\) with three
  distinct realized \(\phi \)-types over \(N\) extending \(\phi
  (x;\gamma )\)).  Thus the only non-isolated \(\phi
  \)-type over \(\mathfrak{C}\) is \(\set{\neg\phi
  (x;c):c\in\mathfrak{C}}\). It is straightforward to check
  that having a unique non-isolated \(\phi \)-type is sufficient to show
  that \(R_\phi (x=x)=\mult_\phi (x=x)=1\).
\end{rpf}

\printbibliography[heading=bibintoc]
\end{document}